\documentstyle{article}

\def\cirk{\,{\raisebox{.3ex}{\tiny $\circ$}}\,}
\def\df{\mbox{\scriptsize{\it df}}}
\def\mj{\mbox{\bf 1}}
\def\prop#1#2{\vspace{2ex} \noindent{\sc #1.} {\it #2} \par \vspace{2ex}}

\begin{document}

\title{Equality of Proofs for Linear Equality}
\author{{\sc Kosta Do\v sen} and {\sc Zoran Petri\' c}\\[0.5cm]
{\small Mathematical Institute} \\[-.5ex]
{\small Knez Mihailova 36, p.f. 367} \\[-.5ex]
{\small 11001 Belgrade, Serbia} \\[-.5ex]
{\small email: \{kosta, zpetric\}@mi.sanu.ac.yu}}
\date{}
\maketitle

\begin{abstract}
\noindent This paper is about equality of proofs in which a binary
predicate formalizing properties of equality occurs, besides
conjunction and the constant true proposition. The properties of
equality in question are those of a preordering relation, those of
an equivalence relation, and other properties appropriate for an
equality relation in linear logic. The guiding idea is that
equality of proofs is induced by coherence, understood as the
existence of a faithful functor from a syntactical category into a
category whose arrows correspond to diagrams. Edges in these
diagrams join occurrences of variables that must remain the same
in every generalization of the proof. It is found that assumptions
about equality of proofs for equality are parallel to standard
assumptions about equality of arrows in categories. They reproduce
standard categorial assumptions on a different level. It is also
found that assumptions for a preordering relation involve an
adjoint situation.
\end{abstract}

\vspace{.3cm}

\noindent {\small {\it Mathematics Subject Classification} ({\it
2000}): 03F07, 03F05, 03F52, 03G30, 18C05, 18A40, 18D10}

\vspace{1ex}

\noindent {\small {\it Keywords$\,$}: equality, linear logic,
preordering relation, equivalence relation, congruence relation,
equality of proofs, categorial coherence, generality of proofs,
monoidal categories, adjunction}

\vspace{1ex}

\noindent{\small {\it Acknowledgements}. We would like to thank
Andreas Blass and an anonymous referee for useful suggestions to
improve the exposition. Work on this paper was supported by the
Ministry of Science of Serbia (Grant 144013).}

\section{Introduction}

The purpose of this paper is to investigate in equational logic
the hypothesis that two proofs are equal if and only if they have
the same generality. Two proofs, with the same premises and
conclusions, have the same generality when after diversifying
variables as much as possible without changing the rules of
inference one ends up again with the same premises and
conclusions, up to renaming of variables. This notion of
generality of proof was investigated in \cite{DP03a} (see also
\cite{DP03b}), \cite{DP04} and \cite{DP05}. However, individual
variables (as opposed to propositional variables) and equality
between them are mentioned only briefly, as an example, in
\cite{DP03a} (end of Section~3). The present paper develops more
fully these matters announced in that previous paper.

The results of this paper take the form of coherence theorems,
understood as faithfulness results for functors from syntactically
constructed categories to a category whose arrows correspond to
diagrams. Edges in these diagrams join occurrences of variables
that must remain the same in every generalization of the proof.

The central result of the paper is that assumptions about equality
of proofs for equality induced by generality are parallel to
standard assumptions about equality of arrows in categories. Usual
categorial assumptions are reproduced on a different level.
Reflexivity of equality corresponds to identity arrows, and
transitivity corresponds to categorial composition of arrows.
Equations involving reflexivity and transitivity are parallel to
the categorial assumptions of omission of identity arrows in
composition and associativity of composition. Similar
correspondences hold for other postulates concerning equality.

Another result of the paper is that assumptions for a preordering
relation involve an adjoint situation. This was foreshadowed in
\cite{L70}.

Besides \cite{L70}, another paper about equality of proofs for
equality is \cite{H87}, which investigates the equivalence of
various syntactical formulations in classical predicate logic with
equality. There is a chapter on equational logic based on
fibrations in \cite{J99} (Chapter 3), but, contrary to what we
have in this paper, it is asserted there (p.\ 174) that there are
no different proofs of the same proposition.

We restrict ourselves to equality added to multiplicative
conjunctive propositional linear logic, with the multiplicative
constant true proposition. We eschew going beyond this limited
fragment of equational logic because generality of proofs
involving equality, like that of proofs involving implication,
prevents the arrows corresponding to the structural rules of
contraction and thinning from making natural transformations (cf.\
\cite{DP02}, Section~1, \cite{DP04}, Section 14.3, and
\cite{DP05a}). This requirement of naturality is otherwise quite
natural, and proof-theoretically well motivated (see \cite{DP04},
Chapters 9-11). Contraction and thinning are required if we want
to say that we deal with full equational logic. In their absence,
we cannot pretend to cover more than a relation of equality
appropriate for linear logic (such as the equality relations
investigated in \cite{D96}). This explains the expression
\emph{linear equality} in the title of the paper.

We restrict ourselves in general to a context with minimal
assumptions where our results can be obtained. So we do not assume
the commutativity of multiplicative conjunction if this is not
essential; i.e., we work also in noncommutative linear logic
(which is related to the Lambek calculus; we stay however within
the multiplicative conjunctive fragment of this calculus).

As we indicated above concerning adjunction, our results are about
relations more general than equality relations, such as
preordering and equivalence relations. Most of the paper (Sections
2-7) is about such relations. Only in the last section we indicate
how to deal with further assumptions, such as congruence. Within
our syntactical systems, where the means of expression are
limited, the equivalence relations involved amount however to
equality relations. In the last section our proofs will be less
formal.

We consider binary operational expressions in the last section,
but we do not have anywhere predicate variables. So we cannot say
that we deal yet with full linear equational logic, but only with
fragments of it. In the context where we investigate equality of
proofs involving equality, motivated by generality, we would have
to enter into the question of what is a linear predicate (see
\cite{D96a}), and moreover we would have to restrict ourselves to
the multiplicative fragment of linear logic without propositional
constants (see \cite{DP05}). Equations that can be expected in
that context would be on the lines of \cite{H87} (Sections 2.3-4).
We leave however these extensions of our approach for another
occasion.

\section{The category $M_\leq$}

Let $\cal V$ be a set whose elements, for which we use the letters
$x,y,z,\ldots$, perhaps with indices, are called \emph{variables}.
The cardinality of $\cal V$ is not restricted: $\cal V$ can be
infinite or finite, and even empty. Let words of the form ${x\leq
y}$ or $\top$ be called \emph{atomic formulae}. The set of
\emph{formulae} is defined inductively as follows. Atomic formulae
are formulae, and if $A$ and $B$ are formulae, then ${(A\wedge
B)}$ is a formula. We use $A,B,C,\ldots$ for formulae, and we
omit, as usual, the outermost parentheses of ${(A\wedge B)}$. (We
proceed analogously for other similar expressions later on.)

The objects of the category $M_\leq$ are formulae. To define the
arrows of $M_\leq$, we define first inductively the \emph{arrow
terms} of $M_\leq$ in the following way. We use $f,g,h,\ldots$,
perhaps with indices, for arrow terms. Every arrow term $f$ has a
\emph{type}, which is an ordered pair ${(A,B)}$ of objects of
$M_\leq$; that $f$ is of type ${(A,B)}$ is written ${f\!:A\vdash
B}$.

For all formulae $A$, $B$ and $C$, and for all variables $x$, $y$
and $z$, the following are \emph{primitive arrow terms} of
$M_\leq$:

\begin{tabbing}
\centerline{$\mj_A\!: A\vdash A$,}
\\[1.5ex]
\mbox{\hspace{1em}}\= $b^{\rightarrow}_{A,B,C}\,$\= : \=
$A\wedge(B\wedge C)\vdash (A\wedge B)\wedge C$,\quad \=
$b^{\leftarrow}_{A,B,C}\,$\= : \= $(A\wedge B)\wedge C\vdash
A\wedge (B\wedge C)$,\kill

\> $b^{\rightarrow}_{A,B,C}\,$\> : \> $A\wedge(B\wedge C)\vdash
(A\wedge B)\wedge C$,\> $b^{\leftarrow}_{A,B,C}\,$\> : \>
$(A\wedge B)\wedge C\vdash A\wedge (B\wedge C)$,
\\[.5ex]
\> \>$\delta^{\rightarrow}_A\!\!$\': \> $A\wedge \top\vdash A$,\>
\> $\delta^{\leftarrow}_A\!\!$\': \> $A\vdash A\wedge \top$,
\\[.5ex]
\> \>$\sigma^{\rightarrow}_A\!\!$\' : \> $\top\wedge A\vdash A$,\>
\> $\sigma^{\leftarrow}_A\!\!$\' : \> $A\vdash \top\wedge A$,
\\[1.5ex]
\centerline{$r_x\!:\top\vdash x\leq x$,}
\\[.5ex]
\centerline{$t_{x,y,z}\!:x\leq y\wedge y\leq z\vdash x\leq z$.}
\end{tabbing}

If ${f\!:A\vdash B}$ and ${g\!:C\vdash D}$ are arrow terms of
$M_\leq$, then ${f\wedge g\!:A\wedge C\vdash B\wedge D}$ and
${g\cirk f\!:A\vdash D}$ are arrow terms of $M_\leq$, provided
that for ${g\cirk f}$ we have that $B$ and $C$ are the same
formula. This defines the arrow terms of $M_\leq$.

The arrows of $M_\leq$ are equivalence classes of arrow terms with
respect to the smallest equivalence relation which guarantees that
the following equations are satisfied:

\begin{tabbing}
\quad\emph{categorial equations:}
\\*[1ex]
\hspace{3em}\= (\emph{cat}~1)\quad\= $\mj_B\cirk f=f,\quad
f\cirk\mj_A=f$,\quad for $f\!:A\vdash B$,
\\*[.5ex]
\> (\emph{cat}~2)\> $(h\cirk g)\cirk f=h\cirk(g\cirk f)$,
\\[2ex]
\quad\emph{bifunctoriality equations:}
\\*[1ex]
\> ($\wedge\:$1)\> $\mj_A\wedge\mj_B=\mj_{A\wedge B}$,
\\*[.5ex]
\> ($\wedge\:$2)\> $(g_1\cirk f_1)\wedge(g_2\cirk f_2)=(g_1\wedge
g_2)\cirk(f_1\wedge f_2)$,
\\[2ex]
\quad\emph{naturality equations:}
\\*[1ex]
\quad\quad for $f\!:A\vdash D$, $g\!: B\vdash E$ and $h\!:C\vdash
F$,
\\*[.5ex]
\> ($b$~\emph{nat})\> $((f\wedge g)\wedge h)\cirk
b^{\rightarrow}_{A,B,C}=b^{\rightarrow}_{D,E,F}\cirk(f\wedge(g\wedge
h))$,
\\[.5ex]
\> ($\delta$~\emph{nat})\>
$f\cirk\delta^{\rightarrow}_A=\delta^{\rightarrow}_D\cirk(f\wedge\mj_\top)$,
\\*[.5ex]
\> ($\sigma$~\emph{nat})\>
$f\cirk\sigma^{\rightarrow}_A=\sigma^{\rightarrow}_D\cirk(\mj_\top\wedge
f)$,
\\[2ex]
\quad\emph{specific equations of monoidal categories:}
\\*[1ex]
\> ($bb$)\> $b^{\leftarrow}_{A,B,C}\cirk
b^{\rightarrow}_{A,B,C}=\mj_{A\wedge(B\wedge C)}$,\quad\quad
$b^{\rightarrow}_{A,B,C}\cirk b^{\leftarrow}_{A,B,C}=\mj_{(A\wedge
B)\wedge C}$,
\\*[.5ex]
\> ($b\:$5)\> $b^{\rightarrow}_{A\wedge B,C,D}\cirk
b^{\rightarrow}_{A,B,C\wedge
D}=(b^{\rightarrow}_{A,B,C}\wedge\mj_D)\cirk
b^{\rightarrow}_{A,B\wedge C,D}\cirk(\mj_A\wedge
b^{\rightarrow}_{B,C,D})$,
\\[.7ex]
\> ($\delta\delta$)\>
$\delta^{\leftarrow}_A\cirk\delta^{\rightarrow}_A\;$\=$=\mj_{A\wedge\top}$,\quad\quad
$\delta^{\rightarrow}_A\cirk\delta^{\leftarrow}_A\;$\=$=\mj_A$,
\\*[.5ex]
\> ($\sigma\sigma$)\>
$\sigma^{\leftarrow}_A\cirk\sigma^{\rightarrow}_A$\>$=\mj_{\top\wedge
A}$,\quad\quad
$\sigma^{\rightarrow}_A\cirk\sigma^{\leftarrow}_A$\>$=\mj_A$,
\\[.7ex]
\> ($b\delta\sigma$)\>
$b^{\rightarrow}_{A,\top,C}=(\delta^{\leftarrow}_A\wedge\mj_C)\cirk(\mj_A\wedge\sigma^{\rightarrow}_C)$,
\\[2ex]
\quad\emph{specific equations of $M_\leq$:}
\\*[1ex]
\> ($rt\delta$)\> $t_{x,y,y}\cirk(\mj_{x\leq y}\wedge
r_y)=\delta^{\rightarrow}_{x\leq y}$,
\\*[.5ex]
\> ($rt\sigma$)\> $t_{y,y,x}\cirk(r_y\wedge \mj_{y\leq
x})=\sigma^{\rightarrow}_{y\leq x}$,
\\[.7ex]
\> ($tb$)\> $t_{x,y,u}\cirk(\mj_{x\leq y}\wedge
t_{y,z,u})=t_{x,z,u}\cirk(t_{x,y,z}\wedge\mj_{z\leq u})\cirk
b^{\rightarrow}_{x\leq y,y\leq z,z\leq u}$;

\end{tabbing}

\noindent if ${f_1=g_1}$ and ${f_2=g_2}$, then ${f_1\wedge
f_2=g_1\wedge g_2}$ and ${f_2\cirk f_1=g_2\cirk g_1}$, provided
${f_2\cirk f_1}$ and ${g_2\cirk g_1}$ are defined.

The category $M_\leq$ is a monoidal category in the sense of
\cite{ML71} (Section VII.1, see also \cite{DP04}, Section 4.6).
The equations ($rt\delta$) and $(rt\sigma$) are parallel to the
equations (\emph{cat}~1), and the equation ($tb$) is parallel in
the same manner to the equation (\emph{cat}~2). The equation
($tb$), for instance, says that a composition tied to the arrows
$t$ is associative.

The following instance of ($tb$):

\[
t_{x,x,x}\cirk(\mj_{x\leq x}\wedge
t_{x,x,x})=t_{x,x,x}\cirk(t_{x,x,x}\wedge\mj_{x\leq x})\cirk
b^{\rightarrow}_{x\leq x,x\leq x,x\leq x}
\]

\noindent is analogous to the equation (\emph{\v b\v w}) of
categories with coproducts (see \cite{DP04}, List of Equations),
where ${t_{x,x,x}\!:x\leq x\wedge x\leq x\vdash x\leq x}$
corresponds to the codiagonal arrow ${\mbox{\emph{\v w}}_p\!:p\vee
p\vdash p}$ of these categories. The instances of ($rt\delta$) and
($rt\sigma$) with $x$ and $y$ the same variable correspond in an
analogous manner to the equations (\emph{\v w\v k}) of categories
with coproducts (see \cite{DP04}, ibid.). The arrow
${r\!:\top\vdash x\leq x}$ corresponds here to the arrow
${\check\kappa_p\!:\bot\vdash p}$ of categories with coproducts,
$\bot$ being the initial object.

The arrows $r_x$ and $t_{x,y,z}$ codify of course the reflexivity
and transitivity of a relation corresponding to $\leq$. So we deal
here with a \emph{preordering} relation.

\section{The coherence of $M_\leq$}

For every object $A$ of $M_\leq$ let $GA$ be the number of
occurrences of variables in $A$. (One could modify the category
\emph{Br}, mentioned below, so that its objects are formulae,
instead of finite ordinals. In that case the object $GA$ would be
the formula $A$. We draw diagrams in this paper in that spirit.
The category \emph{Br} abstracts from a formula just the position
of occurrences of variables in a formula, which is the only thing
relevant for drawing diagrams. We do not expect the functor $G$
below to be one-one on objects. It suffices for our purposes that
it be faithful.)

Let us assign the following diagrams to the primitive arrow terms
of $M_\leq$

\begin{center}
\begin{picture}(30,40)
\put(20,10){\line(0,1){20}}

\put(20,7){\makebox(0,0)[t]{$A$}}
\put(20,33){\makebox(0,0)[b]{$A$}}
\put(15,20){\makebox(0,0)[r]{$\mj_A$}}
\end{picture}
\end{center}

\noindent (where the line joining the two $A$'s stands for a
family of parallel lines---one line for each occurrence of a
variable in $A$; for example, for $A$ being ${x\leq y\wedge z\leq
x}$ we have

\begin{center}
\begin{picture}(80,40)
\put(13,10){\line(0,1){19}} \put(31,10){\line(0,1){19}}
\put(48,10){\line(0,1){19}} \put(66,10){\line(0,1){19}}

\put(40,7){\makebox(0,0)[t]{$x\leq y\wedge z\leq x$}}
\put(40,33){\makebox(0,0)[b]{$x\leq y\wedge z\leq x$}}
\end{picture}
\end{center}

\noindent and analogously in other cases below)

\begin{center}
\begin{picture}(200,40)(-2,0)
\put(25,10){\line(0,1){19}} \put(48,10){\line(0,1){19}}
\put(70,10){\line(0,1){19}}

\put(47,7){\makebox(0,0)[t]{$(A\,\wedge\, B)\wedge C$}}
\put(50,33){\makebox(0,0)[b]{$A\wedge(B\,\wedge\, C)$}}
\put(20,20){\makebox(0,0)[r]{$b^{\rightarrow}_{A,B,C}$}}

\put(125,10){\line(0,1){19}} \put(148,10){\line(0,1){19}}
\put(170,10){\line(0,1){19}}

\put(147,33){\makebox(0,0)[b]{$(A\,\wedge\, B)\wedge C$}}
\put(150,7){\makebox(0,0)[t]{$A\wedge(B\,\wedge\, C)$}}
\put(120,20){\makebox(0,0)[r]{$b^{\leftarrow}_{A,B,C}$}}
\end{picture}
\end{center}

\vspace{1ex}

\begin{center}
\begin{picture}(200,40)
\put(50,10){\line(0,1){20}}

\put(50,7){\makebox(0,0)[t]{$A$}}
\put(60,33){\makebox(0,0)[b]{$A\wedge\top$}}
\put(20,20){\makebox(0,0)[r]{$\delta^{\rightarrow}_A$}}

\put(150,10){\line(0,1){20}}

\put(160,7){\makebox(0,0)[t]{$A\wedge\top$}}
\put(150,33){\makebox(0,0)[b]{$A$}}
\put(120,20){\makebox(0,0)[r]{$\delta^{\leftarrow}_A$}}
\end{picture}
\end{center}

\vspace{1ex}

\begin{center}
\begin{picture}(200,40)(20,0)
\put(70,10){\line(0,1){20}}

\put(70,7){\makebox(0,0)[t]{$A$}}
\put(60,33){\makebox(0,0)[b]{$\top\wedge A$}}
\put(40,20){\makebox(0,0)[r]{$\sigma^{\rightarrow}_A$}}

\put(170,10){\line(0,1){20}}

\put(160,7){\makebox(0,0)[t]{$\top\wedge A$}}
\put(170,33){\makebox(0,0)[b]{$A$}}
\put(140,20){\makebox(0,0)[r]{$\sigma^{\leftarrow}_A$}}
\end{picture}
\end{center}

\vspace{1ex}

\begin{center}
\begin{picture}(200,40)
\put(50,10){\oval(20,20)[t]}

\put(50,7){\makebox(0,0)[t]{$x\leq x$}}
\put(50,33){\makebox(0,0)[b]{$\top$}}
\put(20,20){\makebox(0,0)[r]{$r_x$}}

\put(142,10){\line(-1,1){20}} \put(162,10){\line(1,1){20}}
\put(150,29){\oval(22,22)[b]}

\put(152,7){\makebox(0,0)[t]{$x\leq z$}}
\put(152,33){\makebox(0,0)[b]{$x\leq y\:\wedge\; y\,\leq\, z$}}
\put(115,20){\makebox(0,0)[r]{$t_{x,y,z}$}}
\end{picture}
\end{center}

In the diagram for $r_x$ we have a \emph{cap} joining the two
occurrences of $x$ and in the diagram for $t_{x,y,z}$ we have a
\emph{cup} joining the two occurrences of $y$. We use an analogous
terminology in other cases.

These diagrams and the function $G$ on objects serve to define a
functor $G$ from $M_\leq$ to the category \emph{Br} of \cite{DP05}
(Section 2.3). Namely, $G$ maps an arrow of $M_\leq$ to an arrow
of \emph{Br} that corresponds to a diagram. The composition of
$M_\leq$ is mapped to composition in \emph{Br}, which corresponds
to vertical composition of diagrams, while the operation $\wedge$
on the arrows of $M_\leq$ is mapped to the operation of \emph{Br}
that corresponds to horizontal composition of diagrams (see
\cite{DP05}, Section 2.3). We check by induction on the length of
derivation that $G$ is indeed a functor. Here is what we have in
the basis of this induction for the specific equations of
$M_\leq$:

\begin{center}
\begin{picture}(200,80)
\put(-30,75){\makebox(0,0)[r]{($rt\delta$):}}

\put(40,10){\line(-1,1){20}} \put(60,10){\line(1,1){19}}
\put(18,40){\line(0,1){20}} \put(38,40){\line(0,1){19}}
\put(48,29){\oval(20,20)[b]} \put(70,40){\oval(22,22)[t]}

\put(50,7){\makebox(0,0)[t]{$x\leq y$}}
\put(50,33){\makebox(0,0)[b]{$x\leq y\:\wedge\; y\,\leq\, y$}}
\put(45,63){\makebox(0,0)[b]{$x\leq y\:\wedge\;\;\;\;\top$}}
\put(10,18){\makebox(0,0)[r]{$t_{x,y,y}$}}
\put(10,50){\makebox(0,0)[r]{$\mj_{x\leq y}\wedge r_y$}}

\put(140,30){\line(0,1){20}} \put(159,30){\line(0,1){19}}

\put(150,27){\makebox(0,0)[t]{$x\leq y$}}
\put(160,53){\makebox(0,0)[b]{$x\leq y\wedge\top$}}
\put(130,40){\makebox(0,0)[r]{$\delta^{\rightarrow}_{x\leq y}$}}
\end{picture}
\end{center}

\vspace{1ex}

\begin{center}
\begin{picture}(200,80)
\put(-30,75){\makebox(0,0)[r]{($rt\sigma$):}}

\put(40,10){\line(-1,1){20}} \put(60,10){\line(1,1){20}}
\put(59,40){\line(0,1){19}} \put(80,40){\line(0,1){19}}
\put(48,29){\oval(20,20)[b]} \put(28,40){\oval(22,22)[t]}

\put(50,7){\makebox(0,0)[t]{$y\leq x$}}
\put(50,33){\makebox(0,0)[b]{$y\leq y\:\wedge\; y\,\leq\, x$}}
\put(54,63){\makebox(0,0)[b]{$\top\;\;\;\;\wedge\;\,y\,\leq\, x$}}
\put(10,18){\makebox(0,0)[r]{$t_{y,y,x}$}}
\put(10,50){\makebox(0,0)[r]{$r_y\wedge\mj_{y\leq x}$}}

\put(160,30){\line(0,1){19}} \put(179,30){\line(0,1){19}}

\put(170,27){\makebox(0,0)[t]{$y\leq x$}}
\put(160,53){\makebox(0,0)[b]{$\top\wedge y\leq x$}}
\put(130,40){\makebox(0,0)[r]{$\sigma^{\rightarrow}_{y\leq x}$}}
\end{picture}
\end{center}

\vspace{1ex}

\begin{center}
\begin{picture}(300,120)
\put(20,115){\makebox(0,0)[r]{($tb$):}}

\put(74,20){\line(-4,5){15}} \put(95,20){\line(4,5){15}}
\put(83,39){\oval(20,20)[b]}

\put(51,50){\line(-4,5){15}} \put(71,50){\line(-4,5){15}}
\put(92,50){\line(-4,5){15}} \put(117,50){\line(4,5){15}}
\put(102,69){\oval(20,20)[b]}

\put(85,17){\makebox(0,0)[t]{$x\leq u$}}
\put(85,43){\makebox(0,0)[b]{$x\leq y\:\wedge\; y\,\leq\, u$}}
\put(85,73){\makebox(0,0)[b]{$x\leq y\wedge(y\leq z\;\wedge\;
z\leq u)$}} \put(38,27){\makebox(0,0)[r]{$t_{x,y,u}$}}
\put(38,57){\makebox(0,0)[r]{$\mj_{x\leq y}\wedge t_{y,z,u}$}}

\put(249,10){\line(-4,5){15}} \put(270,10){\line(4,5){15}}
\put(258,29){\oval(20,20)[b]}

\put(226,40){\line(-4,5){15}} \put(250,40){\line(4,5){15}}
\put(272,40){\line(4,5){15}} \put(292,40){\line(4,5){15}}
\put(239,59){\oval(20,20)[b]}

\put(211,73){\line(0,1){20}} \put(230,73){\line(0,1){20}}
\put(249,73){\line(0,1){20}} \put(268,73){\line(0,1){20}}
\put(291,73){\line(0,1){20}} \put(310,73){\line(0,1){20}}

\put(260,7){\makebox(0,0)[t]{$x\leq u$}}
\put(260,33){\makebox(0,0)[b]{$x\leq z\:\wedge\; z\,\leq\, u$}}
\put(260,63){\makebox(0,0)[b]{$(x\leq y\,\wedge\, y\leq
z)\,\wedge\, z\leq u\:$}}

\put(263,97){\makebox(0,0)[b]{$x\leq y\wedge (y\leq z\;\wedge\;
z\leq u)$}}

\put(207,17){\makebox(0,0)[r]{$t_{x,z,u}$}}
\put(207,47){\makebox(0,0)[r]{$t_{x,y,z}\wedge\mj_{z\leq u}$}}
\put(207,87){\makebox(0,0)[r]{$b^{\rightarrow}_{x\leq y,y\leq
z,z\leq u}$}}

\end{picture}
\end{center}

\noindent (Note that there are no cups and caps in the diagrams
corresponding to the left-hand sides of ($rt\delta$) and
($rt\sigma$); they were abolished by composing.)

So if $f=g$ in $M_\leq$, then ${Gf=Gg}$. Our purpose is to show
also the converse for $f$ and $g$ of the same type; we show,
namely, that $G$ is a faithful functor from $M_\leq$ to \emph{Br}.
We call this faithfulness result the \emph{coherence} of $M_\leq$.

An arrow term of $M_\leq$ of the form ${f_n \cirk \ldots \cirk
f_1}$, where $n\geq 1$, with parentheses tied to $\cirk$
associated arbitrarily, such that for every ${i\in\{1,\ldots,n\}}$
we have that $f_i$ is without $\cirk$ is called \emph{factorized}.
In a factorized arrow term ${f_n\cirk\ldots\cirk f_1}$ the arrow
terms $f_i$ are called \emph{factors}.

If $\beta$ is a primitive arrow term of $M_\leq$ which is not of
the form $\mj_B$, then $\beta$-\emph{terms} are defined
inductively as follows: $\beta$ is a $\beta$-term; if $f$ is a
$\beta$-term, then for every object $A$ of $M_\leq$ we have that
${\mj_A\wedge f}$ and ${f\wedge \mj_A}$ are $\beta$-terms. In a
$\beta$-term the subterm $\beta$ is called the \emph{head} of this
$\beta$-term. For example, the head of the $t_{x,y,z}$-term
${(\mj_{u\leq v}\wedge t_{x,y,z})\wedge\mj_{x\leq u}}$ is
$t_{x,y,z}$.

We define $\mj$-\emph{terms} as $\beta$-terms by replacing $\beta$
in the definition above by $\mj_B$.

A factor that is a $\beta$-term for some $\beta$ is called a
\emph{headed} factor. A factorized arrow term is called
\emph{headed} when each of its factors is either headed or a
$\mj$-term. A headed arrow term ${f_n\cirk\ldots\cirk f_1}$ is
called \emph{developed} when $f_1$ is a $\mj$-term and if $n>1$,
then every factor of ${f_n\cirk\ldots\cirk f_2}$ is headed.
Analogous definitions of factorized arrow term, factor,
$\beta$-term, head, headed factor, headed factorized arrow term
and developed arrow term can be given later for categories other
than $M_\leq$, and we will not dwell on these definitions any
more.

By using the categorial equations (\emph{cat}~1) and
(\emph{cat}~2) and the bifunctoriality equations we can easily
prove by induction on the length of $f$ the following lemma for
$M_\leq$.

\prop{Development Lemma}{For every arrow term $f$ there is a
developed arrow term $f'$ such that ${f=f'}$.}

An $r$-\emph{less} arrow term of $M_\leq$ is an arrow term of
$M_\leq$ in which $r_x$ does not occur for any $x$. A headed
factorized arrow term each of whose factors is an $r_x$-term for
some $x$ or a $\mj$-term is called an $r$-\emph{factorized} arrow
term. We can easily prove the following lemma by applying the
Development Lemma and the equations ($rt\delta$) and ($rt\sigma$),
besides bifunctoriality, naturality and other obvious equations.

\prop{$r$-Normality Lemma}{For every arrow term $f$ of $M_\leq$
there is a headed factorized arrow term of $M_\leq$ of the form
${f_r\cirk f'}$ such that $f$ is a developed $r$-less arrow term
and $f_r$ is an $r$-factorized arrow term and ${f=f_r\cirk f'}$ in
$M_\leq$.}

\noindent The arrow term ${f_r\cirk f'}$ of this lemma is called
the $r$-\emph{normal form} of $f$.

Suppose ${f,g\!:A\vdash B}$ are arrow terms of $M_\leq$ such that
${Gf=Gg}$. Let ${f_r\cirk f'}$ and ${g_r\cirk g'}$ be the
$r$-normal forms of $f$ and $g$ respectively. Then $G(f_r\cirk
f')=G(g_r\cirk g')$, and there is a bijection between the caps of
${G(f_r\cirk f')}$ and ${G(g_r\cirk g')}$. Moreover, there is a
bijection between the caps in ${G(f_r\cirk f')}$ and the
$r_x$-factors of $f_r$, and analogously for ${G(g_r\cirk g')}$ and
$g_r$.

By using the bifunctoriality equations we can achieve that $f_r$
and $g_r$, which are $r$-factorized, are the same arrow term $h$.
Since

\[
\begin{array}{l}
G(h\cirk f')=Gh\cirk Gf',
\\[.5ex]
G(h\cirk g')=Gh\cirk Gg',
\\[.5ex]
G(h\cirk f')=G(h\cirk g'),
\end{array}
\]

\noindent and $Gh$ has no cups, it is easy to conclude that

\[
Gf'=Gg'.
\]

\noindent (In the category \emph{Br} the arrow $Gh$, which has no
cups, has a left inverse, which is its image in a horizontal
mirror.)

There are no caps in $Gf'$ and $Gg'$, and there is a bijection
between the cups of $Gf'$ and the $t_{x,y,z}$-factors of $f'$, and
analogously for $Gg'$ and $g'$. A $t_{x,y,z}$-factor of $f'$ and a
$t_{u,y,v}$-factor of $g'$ that correspond to each other according
to these bijections are called \emph{coupled}. Suppose $Gf'$ and
$Gg'$ have at least one cup. Then by using the equation ($tb$),
besides the bifunctoriality, naturality and other obvious
equations, we obtain ${f'=h_1\cirk f''}$ and ${g'=h_2\cirk g''}$
where $f''$ and $g''$ are developed arrow terms, while $h_1$ is a
$t_{x,y,z}$-factor coupled with the $t_{u,y,v}$-factor $h_2$.

It is impossible that $z$ coincides with $v$ while $x$ differs
from $u$. Otherwise, the targets of $h_1$ and $h_2$ would differ.
If $x$ differs from $u$, and $z$ differs from $v$, then in the
source of $f''$ and $g''$ we would have variables occurring in the
following order:

\[
{x\ldots z\ldots u\ldots v}\quad \mbox{\rm{ or }} \quad{u\ldots
v\ldots x\ldots z},
\]

\noindent with a cup between $x$ and $z$ in ${G(h_1\cirk f'')}$
corresponding to $h_1$ and a cup between $u$ and $v$ in
${G(h_2\cirk g'')}$ corresponding to $h_2$. This is impossible
because $h_1$ and $h_2$ are coupled. So $t_{x,y,z}$ coincides with
$t_{u,y,v}$, and by using perhaps the bifunctoriality equation
($\wedge\:$1) we can achieve that $h_1$ and $h_2$ are the same
arrow term. From ${G(h\cirk f'')=G(h\cirk g'')}$ we conclude that
${Gf''=Gg''}$.

Then we proceed by induction on the number of cups in $Gf'$, which
is equal to $Gg'$, to show that ${f'=g'}$. In the basis of this
induction we rely on Mac Lane's monoidal coherence (see
\cite{ML63} and \cite{ML71}, Section VII.2, or \cite{DP04},
Section 4.6). From that it follows that if for ${f,g\!:A\vdash B}$
arrow terms of $M_\leq$ we have ${Gf=Gg}$, then ${f=g}$ in
$M_\leq$. This proves the coherence of $M_\leq$.

Note that by omitting from the proof above the part involving
$r_x$ we would obtain an analogous coherence result for a category
defined like $M_\leq$ save that it lacks the arrows $r_x$ and the
specific equations ($rt\delta$) and ($rt\sigma$). We can also
obtain coherence for a category defined like $M_\leq$ but lacking
the arrows $t_{x,y,z}$ and all the specific equations of $M_\leq$.

\section{The coherence of $S_\leq$}

The category $S_\leq$ is defined like $M_\leq$ save that we have
for all formulae $A$ and $B$ the additional primitive arrow term

\[
c_{A,B}\!:A\wedge B\vdash B\wedge A,
\]

\noindent which is subject to the following additional equations:

\begin{tabbing}
\quad\emph{naturality equation:}
\\*[1ex]
\hspace{3em}\= ($c$~\emph{nat})\quad\= $(g\wedge f)\cirk
c_{A,B}=c_{D,E}\cirk (f\wedge g)$,
\\[2ex]
\quad\emph{specific equations of symmetric monoidal categories:}
\\*[1ex]
\> ($cc$)\> $c_{B,A}\cirk c_{A,B}=\mj_{A\wedge B}$,
\\*[.5ex]
\> ($bc$)\> $c_{A,B\wedge
C}=b^{\rightarrow}_{B,C,A}\cirk(\mj_B\wedge c_{A,C})\cirk
b^{\leftarrow}_{B,A,C}\cirk(c_{A,B}\wedge\mj_C)\cirk
b^{\rightarrow}_{A,B,C}$.
\end{tabbing}

The category $S_\leq$ is a symmetric monoidal category (see
\cite{ML71}, Section VII.7), in which $\leq$ corresponds to a
preordering relation.

We define the functor $G$ from $S_\leq$ to \emph{Br} by extending
the definition of $G$ from $M_\leq$ to \emph{Br} with a clause
corresponding to the following diagram:

\begin{center}
\begin{picture}(40,40)
\put(10,10){\line(1,1){20}} \put(30,10){\line(-1,1){20}}

\put(20,7){\makebox(0,0)[t]{$B\wedge A$}}
\put(20,33){\makebox(0,0)[b]{$A\,\wedge\, B$}}
\put(0,20){\makebox(0,0)[r]{$c_{A,B}$}}
\end{picture}
\end{center}

That $G$ is indeed a functor follows from well-known facts about
symmetric monoidal categories (which were established in
\cite{ML63}; see also \cite{DP04}, Chapter~5). We can prove
coherence for $S_\leq$ with respect to $G$, i.e., we can prove
that $G$ is faithful, by imitating the proof of coherence for
$M_\leq$ in the preceding section. The only difference is that we
appeal to symmetric monoidal coherence (see \cite{ML63},
Section~5, second edition of \cite{ML71}, Section XI.1, and
\cite{DP04}, Section 5.3) where we appealed before to monoidal
coherence, and we replace the proof that $h_1$ and $h_2$ can be
taken to be the same arrow term by the following alternative
proof. An analogous proof could have already been used in the
preceding section, but there, in the absence of $c_{A,B}$, we also
had a slightly simpler argument.

For ${f\!:A\vdash B}$ an $r$-less arrow term of $S_\leq$, we say
that a set $U$ of occurrences of variables in $A$ is
$f$-\emph{closed} when the following implication holds: if either
$u\leq u'$ is a subformula of $A$ or $u$ and $u'$ are connected by
a cup of $Gf$, and one of $u$ and $u'$ is in $U$, then the other
is in $U$ too. It is easy to verify by induction on the complexity
of $f$ that
\begin{itemize}
\item[$(\ast)$] for every $f$-closed set $U$ and for every atomic
subformula $x\leq y$ of $B$, a member of $U$ is connected by $Gf$
to $x$ if and only if a member of $U$ is connected by $Gf$ to $y$.
\end{itemize}

This holds in particular for $f$-closed sets generated by a single
occurrence of a variable in $A$. We call such $f$-closed sets
\emph{maximal sequences}. It is easy to see that a maximal
sequence is a set $\{u_1,u_2,\ldots,u_{2n-1},u_{2n}\}$ of
occurrences of variables in $A$ for ${n\geq 1}$ such that
${u_{2i-1}\leq u_{2i}}$ is a subformula of $A$, for ${1\leq i\leq
n}$, while $u_{2j}$ and $u_{2j+1}$, for ${1\leq j\leq n-1}$, are
connected by a cup of $Gf$. Note that it follows from $(\ast)$
that for a maximal sequence there must exist an atomic subformula
${x\leq y}$ of $B$ such that $u_1$ in $A$ is connected by $Gf$ to
the occurrence of $x$ in ${x\leq y}$ in $B$ and $u_{2n}$ in $A$ is
connected by $Gf$ to the occurrence of $y$ in ${x\leq y}$ in $B$.

Suppose, as in the preceding section, that for $r$-less arrow
terms ${f',g'\!:A\vdash B}$ of $S_\leq$, we have ${Gf'=Gg'}$.
Again, by using ($tb$), besides bifunctoriality, naturality and
other obvious equations, we obtain ${f'=h_1\cirk f''}$ and
${g'=h_2\cirk g''}$ where $h_1$ is a $t_{x,y,z}$-factor coupled
with the $t_{u,y,v}$-factor $h_2$. Then for the same maximal
sequence $u_1,\ldots,u_{2n}$ in $A$ we have that $u_1$ is
connected by $Gf'$ to the $x$ of $t_{x,y,z}$ and by $Gg'$ to the
$u$ of $t_{u,y,v}$. Analogously, we have that $u_{2n}$ is
connected by $Gf'$ to the $z$ of $t_{x,y,z}$ and by $Gg'$ to the
$v$ of $t_{u,y,v}$. This means that $t_{x,y,z}$ coincides with
$t_{u,y,v}$, and by using perhaps ($\wedge\:$1) we can achieve
that $h_1$ and $h_2$ are the same arrow term.

\section{The coherence of $M_\equiv$}

The category $M_\equiv$ is defined like $M_\leq$ save that $\leq$
is replaced everywhere by $\equiv$, and we have for all variables
$x$ and $y$ the additional primitive arrow term

\[
s_{x,y}\!:x\equiv y\vdash y\equiv x,
\]

\noindent which is subject to the following additional equations:

\begin{tabbing}
\hspace{3em}\= ($ss$)\quad\= $s_{y,x}\cirk s_{x,y}=\mj_{x\equiv
y}$,
\\*[.5ex]
\> ($rs$)\> $s_{x,x}\cirk r_x=r_x$.
\end{tabbing}

\noindent The category $M_\equiv$ is a monoidal category in which
$\equiv$ corresponds to an equivalence relation. Since the means
of expression of $M_\equiv$ are limited, this equivalence relation
is an equality relation.

As the specific equations of $M_\leq$ are parallel to the
categorial equations (\emph{cat}~1) and (\emph{cat}~2) (see
Section~2), so the equations ($ss$) and ($rs$) are parallel to
equations of groupoids, i.e.\ categories where every arrow $f$ has
an inverse $f^{-1}$. The equation ($ss$) corresponds to
${(f^{-1})^{-1}=f}$ and ($rs$) corresponds to
${\mj^{-1}_A=\mj_A}$.

We define the functor $G$ from $M_\equiv$ to \emph{Br} by
extending $G$ from $M_\leq$ to \emph{Br} with a clause
corresponding to the following diagram:

\begin{center}
\begin{picture}(40,40)
\put(10,10){\line(1,1){19}} \put(30,10){\line(-1,1){20}}

\put(20,7){\makebox(0,0)[t]{$y\,\equiv\, x$}}
\put(20,33){\makebox(0,0)[b]{$x\,\equiv\, y$}}
\put(0,20){\makebox(0,0)[r]{$s_{x,y}$}}
\end{picture}
\end{center}

\noindent Since for the equation $(ss)$ we have

\begin{center}
\begin{picture}(140,80)
\put(10,10){\line(1,1){19}} \put(30,10){\line(-1,1){20}}
\put(10,40){\line(1,1){19}} \put(30,40){\line(-1,1){20}}

\put(20,7){\makebox(0,0)[t]{$x\,\equiv\, y$}}
\put(20,33){\makebox(0,0)[b]{$y\,\equiv\, x$}}
\put(20,68){\makebox(0,0)[t]{$x\,\equiv\, y$}}
\put(0,20){\makebox(0,0)[r]{$s_{y,x}$}}
\put(0,50){\makebox(0,0)[r]{$s_{x,y}$}}

\put(109,30){\line(0,1){20}} \put(130,30){\line(0,1){19}}

\put(120,27){\makebox(0,0)[t]{$x\,\equiv\, y$}}
\put(120,53){\makebox(0,0)[b]{$x\,\equiv\, y$}}
\put(100,40){\makebox(0,0)[r]{$\mj_{x\equiv y}$}}
\end{picture}
\end{center}

\noindent and for $(rs)$ we have

\vspace{-1ex}

\begin{center}
\begin{picture}(140,80)
\put(10,10){\line(1,1){19}} \put(30,10){\line(-1,1){20}}
\put(20,40){\oval(22,22)[t]}

\put(20,7){\makebox(0,0)[t]{$x\,\equiv\, x$}}
\put(20,33){\makebox(0,0)[b]{$x\,\equiv\, x$}}
\put(20,68){\makebox(0,0)[t]{$\top$}}
\put(0,20){\makebox(0,0)[r]{$s_{x,x}$}}
\put(0,50){\makebox(0,0)[r]{$r_x$}}

\put(120,30){\oval(22,22)[t]}

\put(120,27){\makebox(0,0)[t]{$x\,\equiv\, x$}}
\put(120,53){\makebox(0,0)[b]{$\top$}}
\put(100,40){\makebox(0,0)[r]{$r_x$}}
\end{picture}
\end{center}

\noindent we can conclude that $G$ is indeed a functor from
$M_\equiv$ to \emph{Br}. In the remainder of this section we prove
the faithfulness of $G$; namely, the coherence of $M_\equiv$.

This proof is more complex than the other proofs of coherence in
this paper. It involves a number of details. We will mention most
of them, but not all, in order not to prolong the exposition
excessively. We do not consider $M_\equiv$ as the most significant
category of this paper. (We find $S_\equiv$ of the next section
more important, and for it coherence is proved more easily.)

We can easily prove an analogue of the $r$-Normality Lemma of
Section~3 for $M_\equiv$. To prove this analogue we apply also the
equation ($rs$).

Then it is enough to prove coherence for $r$-less arrow terms of
$M_\equiv$ to obtain coherence for the whole of $M_\equiv$. (For
$r$-factorized arrow terms we proceed as in Section~3.)

An arrow term of $M_\equiv$ is called
$\delta$-$\sigma$-\emph{less} when $\delta^{\rightarrow}_A$,
$\delta^{\leftarrow}_A$, $\sigma^{\rightarrow}_A$ or
$\sigma^{\leftarrow}_A$ does not occur in it for any $A$. We can
establish the following.

\prop{$\delta$-$\sigma$-Normality Lemma}{For every $r$-less arrow
term ${f\!:A\vdash B}$ of $M_\equiv$ such that $\top$ does not
occur in ${A\vdash B}$ or both $A$ and $B$ are $\top$, there is a
$\delta$-$\sigma$-less arrow term ${f'\!:A\vdash B}$ such that
${f=f'}$ in $M_\equiv$.}

\noindent For the proof of this lemma we rely on bifunctoriality
and naturality equations, and on monoidal coherence. Intuitively,
we push every $\delta^{\rightarrow}_C$-factor in a headed
factorized arrow term towards the right (or
$\delta^{\leftarrow}_C$-factor towards the left), where it or its
descendant will disappear in virtue of the equations
($\delta\delta$) or ($\sigma\sigma$).

For every formula $A$ we define a formula $A^\dagger$ in which
$\top$ does not occur, or which is $\top$, in the following
inductive manner: if $A$ is atomic, then $A^\dagger$ is $A$, and
if $A$ is $B\wedge C$, then $(B\wedge C)^\dagger$ is either
$B^\dagger\wedge C^\dagger$ when neither $B^\dagger$ nor
$C^\dagger$ is $\top$, or $B^\dagger$ when $C^\dagger$ is $\top$,
or $C^\dagger$ when $B^\dagger$ is $\top$. It is clear that there
is an isomorphism $\varphi_A\!:A\vdash A^\dagger$. For every arrow
${f\!:A\vdash B}$ of $M_\equiv$, let ${f^\dagger\!:A^\dagger\vdash
B^\dagger}$ be the arrow ${\varphi_B\cirk f\cirk\varphi^{-1}_A}$.
We have that ${Gf=Gf^\dagger}$, and $f=g$ if and only if
${f^\dagger=g^\dagger}$.

A type ${A\vdash B}$ is \emph{diversified} when every variable in
it occurs exactly twice (once in $A$ and once in $B$, or twice in
$A$, or twice in $B$). An arrow term whose type is diversified is
also called \emph{diversified}.

For every arrow term ${f\!:A\vdash B}$ of $M_\equiv$ there is a
diversified arrow term ${f'\!:A'\vdash B'}$ of $M_\equiv$ such
that $f$ is obtained from $f'$ by substitution in the variables of
$f'$. (Here variables are uniformly replaced by variables.) This
is clear from $Gf$, which dictates how the diversification is to
be achieved. If ${f,g\!:A\vdash B}$ are diversified arrow terms of
$M_\equiv$, then ${Gf=Gg}$. We also have that for ${f,g\!:A\vdash
B}$, there are diversified arrow terms ${f',g'\!:A'\vdash B'}$
such that $f$ and $g$ are substitution instances of $f'$ and $g'$
respectively if and only if ${Gf=Gg}$.

For a headed factorized arrow term ${f_n\cirk\ldots\cirk f_1}$ of
$M_\equiv$, whose factors are $f_1,\ldots,f_n$, we have that
$Gf_i$, for ${1\leq i\leq n}$, contains a crossing if and only if
$f_i$ is an $s_{x,y}$-factor.

For ${f\!:A\vdash B}$ an arrow term of $M_\equiv$ we say that a
cup in the diagram corresponding to $Gf$ \emph{covers} an
occurrence of $\wedge$ in $A$ when the ends of this cup are on
different sides of this occurrence of $\wedge$. For example, in

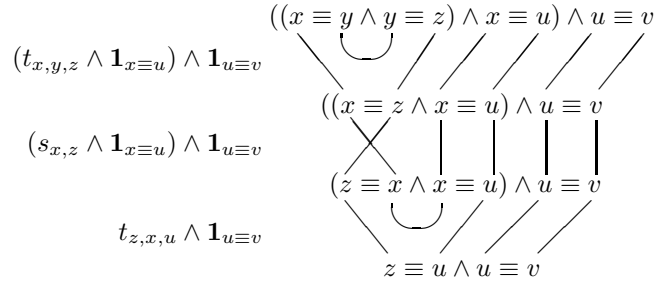
\begin{center}
\begin{picture}(250,110)
\put(122,10){\line(-4,5){16}} \put(142,10){\line(4,5){16}}
\put(159,10){\line(1,1){20}} \put(179,10){\line(1,1){20}}
\put(133,29){\oval(19,19)[b]}

\put(106,40){\line(4,5){17}} \put(125,40){\line(-4,5){17}}
\put(142,40){\line(0,1){20}} \put(162,40){\line(0,1){20}}
\put(182,40){\line(0,1){20}} \put(201,40){\line(0,1){20}}

\put(105,72){\line(-4,5){17}} \put(126,72){\line(2,3){14}}
\put(142,72){\line(4,5){17}} \put(162,72){\line(4,5){17}}
\put(182,72){\line(4,5){17}} \put(202,72){\line(4,5){17}}
\put(114,93){\oval(19,19)[b]}

\put(150,7){\makebox(0,0)[t]{$z\equiv u \wedge u\equiv v$}}
\put(150,33){\makebox(0,0)[b]{$\;(z\equiv x \wedge x\equiv
u)\wedge u\equiv v$}} \put(150,63){\makebox(0,0)[b]{($(x\equiv
z\wedge x\equiv u)\wedge u\equiv v$}}

\put(150,97){\makebox(0,0)[b]{$((x\equiv y \wedge y\equiv z)\wedge
x\equiv u)\wedge u\equiv v$}}

\put(75,17){\makebox(0,0)[r]{$t_{z,x,u}\wedge\mj_{u\equiv v}$}}
\put(75,51){\makebox(0,0)[r]{$(s_{x,z}\wedge \mj_{x\equiv
u})\wedge\mj_{u\equiv v}$}}
\put(75,83){\makebox(0,0)[r]{$(t_{x,y,z}\wedge\mj_{x\equiv
u})\wedge\mj_{u\equiv v}$}}

\end{picture}
\end{center}

\vspace{.75ex}

\noindent the $y$-cup covers only the leftmost occurrence of
$\wedge$ in $((x\equiv y\wedge y\equiv z)\wedge x\equiv u)\wedge
u\equiv v$, and the $x$-cup covers the leftmost and middle
occurrence of $\wedge$. The rightmost occurrence of $\wedge$ is
\emph{uncovered}; i.e., it is not covered by any cup.

Suppose now that for the headed factorized arrow term
${f\!:A\vdash B}$ of $M_\equiv$ we have that  it is $r$-less,
$\delta$-$\sigma$-less and diversified. Then there is an obvious
one-to-one correspondence between occurrences of $\wedge$ in $B$
and uncovered occurrences of $\wedge$ in $A$. There is also an
obvious one-to-one correspondence between the following sets:

\vspace{1.5ex}

the set of $t_{x,y,z}$-factors of $f$,

the set of cups in $Gf$,

the set of variables occurring in $A$ and not in $B$,

the set of occurrences of $\wedge$ in $A$ covered by a cup of
$Gf$.

\vspace{1.5ex}

\noindent Note that all of these one-to-one correspondences that
do not involve the first of these four sets do not depend on the
arrow term $f$, but only on $Gf$.

An arrow term of $M_\equiv$ is called $s$-\emph{normal} when for
every pair of variables ${(x,y)}$ there is at most one occurrence
of $s$ in this arrow term with the indices $_{x,y}$ or $_{y,x}$.

We can easily verify the following.

\prop{$s$-Normality Lemma}{For every diversified arrow term $f$ of
$M_\equiv$ there is a developed $s$-normal arrow term $f'$ of
$M_\equiv$ such that ${f=f'}$ in $M_\equiv$. If $f$ is $r$-less,
then $f'$ is $r$-less too.}

\noindent This holds because, in a diversified developed arrow
term, between two factors whose heads are $s_{x,y}$ or $s_{y,x}$
there can be no factor whose head is $t_{z,x,u}$ or $t_{z,y,u}$,
which would be the only obstacle to bringing the two factors
together, where they get cancelled.

In virtue of all that we have above it is enough to establish the
following in order to prove coherence for $M_\equiv$.

\prop{Auxiliary Lemma}{Suppose $f$ and $g$ are developed,
$r$-less, $\delta$-$\sigma$-less, diversified and $s$-normal arrow
terms of $M_\equiv$ of the same type. Then ${f=g}$ in $M_\equiv$.}

\noindent{\emph{Proof.  }} We proceed by induction on the number
$n$ of $s_{x,y}$-factors and $t_{z,u,v}$-factors in
${f,g\!:A\vdash B}$. This number must be the same in $f$ and $g$
because they are diversified and $s$-normal. (Note that
${Gf=Gg}$.) If ${n=0}$, then we apply monoidal coherence. If
${n>0}$, then there is in $B$ an atomic subformula ${x\equiv y}$
such that either (1) $y$ is in $A$ on the left-hand-side of $x$,
or (2) $x$ is in $A$ on the left-hand side of $y$ and ${x\equiv
y}$ is not a subformula of $A$.

In case (1) we have that

\[
f=h\cirk f'\quad\mbox{\rm and}\quad g=h\cirk g'
\]

\noindent for an $s_{y,x}$-factor $h$, and we may apply the
induction hypothesis to $f'$ and $g'$.

In case (2), we have that

\[
f=h\cirk f'
\]

\noindent for a $t_{x,z,y}$-factor $h$. There must be a
$t_{u,z,v}$-factor $h'$ in $g$. Note that the occurrence of
$\wedge$ in $A$ corresponding to $h$ is covered just by the cup of
$Gf$ corresponding to $h$. This cup in $Gg$, which is equal to
$Gf$, corresponds to $h'$ in $g$. The arrow term $g$ is of the
form

\[
g_m\cirk\ldots\cirk g_1\cirk h'\cirk g'
\]

\noindent for ${m\geq 0}$ (if ${m=0}$, then we have just ${h'\cirk
g'}$); here  ${g_1,\ldots ,g_m}$ are factors. We proceed by
induction on $m$ to show that $g$ is equal to ${h''\cirk g''}$ for
a $t_{x',z,y'}$-factor $h''$, which must be the same as $h$, for
reasons given in the preceding section. (We only replace $\leq$ by
$\equiv$; moreover, ``${u_{2i-1}\leq u_{2i}}$'' is replaced by
``${u_{2i-1}\equiv u_{2i}}$ \emph{or} ${u_{2i}\equiv u_{2i-1}}$''
and ``${x\leq y}$'' is replaced by ``${x\equiv y}$ \emph{or}
${y\equiv x}$''.) Note that there can be no factor in
$g_m\cirk\ldots\cirk g_1$ whose head is $s_{u,v}$, because, as we
said above, the occurrence of $\wedge$ in $A$ corresponding to $h$
is covered just by the cup of $Gf$ corresponding to $h$. \hfill
$\dashv$

\section{The coherence of $S_\equiv$}

The category $S_\equiv$ is defined like $S_\leq$ save that, as
when obtaining $M_\equiv$ out of $M_\leq$, the symbol $\leq$ is
replaced everywhere by $\equiv$ and we have the additional
primitive arrow terms $s_{x,y}$ subject to the equation ($ss$) and
($rs$) of the preceding section, and the additional equation

\begin{tabbing}
\hspace{3em}\= ($ts$)\quad\= $s_{x,z}\cirk
t_{x,y,z}=t_{z,y,x}\cirk(s_{y,z}\wedge s_{x,y})\cirk c_{x\equiv
y,y\equiv z}$.
\end{tabbing}

\noindent This equation is parallel to the following equation of
groupoids (cf.\ Section~5): ${(g\cirk f)^{-1}=f^{-1}\cirk
g^{-1}}$. It is analogous also to the equation (\emph{\v c\v w})
of categories with coproducts (see \cite{DP04}, List of Equations,
and cf.\ the end of Section~2 above).

Note that in the presence of ($ts$) we can derive ($rt\delta$)
from ($rt\sigma$), or vice versa. Here is a derivation of
($rt\delta$) from ($rt\sigma$):

\begin{tabbing}
$t_{x,y,y}\cirk(\mj_{x\equiv y}\wedge r_y)$ \= $=s_{y,x}\cirk
t_{y,y,x}\cirk c_{y\equiv x,y\equiv y}\cirk(s_{x,y}\wedge
s_{y,y})\cirk(\mj_{x\equiv y}\wedge r_y)$, with ($ts$),
\\[1ex]
\> $=s_{y,x}\cirk t_{y,y,x}\cirk(r_y\wedge s_{x,y})\cirk
c_{x\equiv y,\top}$,\quad with ($rs$) and ($c$~\emph{nat}),
\\[1ex]
\> $=s_{y,x}\cirk\sigma^{\rightarrow}_{y\equiv
x}\cirk(\mj_\top\wedge s_{x,y})\cirk c_{x\equiv y,\top}$,\quad
with ($rt\sigma$),
\\[1ex]
\> $=\delta^{\rightarrow}_{x\equiv y}$,\quad with
($\sigma$~\emph{nat}), ($ss$) and monoidal coherence.
\end{tabbing}

The category $S_\equiv$ is a symmetric monoidal category in which
$\equiv$ corresponds to an equivalence relation. Since the means
of expression of $S_\equiv$ are limited, this equivalence relation
is an equality relation.

We define the functor $G$ from $S_\equiv$ to \emph{Br} by
combining what we had for $G$ from $M_\equiv$ to \emph{Br} and for
$G$ from $S_\leq$ to \emph{Br}. Since for the equation $(ts)$ we
have

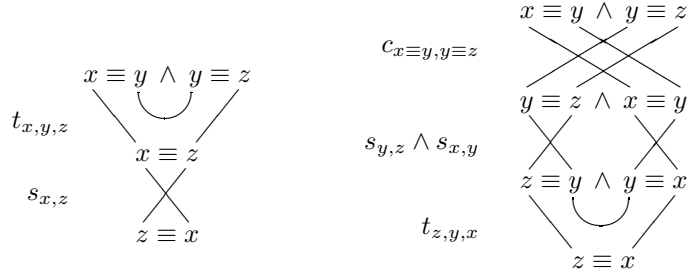
\begin{center}
\begin{picture}(300,110)
\put(66,20){\line(4,5){16}} \put(83,20){\line(-4,5){16}}

\put(63,50){\line(-4,5){16}} \put(86,50){\line(4,5){16}}
\put(74,69){\oval(20,20)[b]}

\put(75,17){\makebox(0,0)[t]{$z\equiv x$}}
\put(75,43){\makebox(0,0)[b]{$x\equiv z$}}
\put(75,73){\makebox(0,0)[b]{$x\equiv y\: \wedge\: y\equiv z$}}
\put(38,29){\makebox(0,0)[r]{$s_{x,z}$}}
\put(38,57){\makebox(0,0)[r]{$t_{x,y,z}$}}

\put(228,10){\line(-4,5){16}} \put(251,10){\line(4,5){16}}
\put(239,29){\oval(21,21)[b]}

\put(212,40){\line(4,5){16}} \put(228,40){\line(-4,5){16}}
\put(252,40){\line(4,5){15}} \put(268,40){\line(-4,5){16}}

\put(211,71){\line(5,3){38}} \put(230,71){\line(5,3){38}}
\put(250,71){\line(-5,3){38}} \put(269,71){\line(-5,3){38}}

\put(240,7){\makebox(0,0)[t]{$z\equiv x$}}
\put(240,33){\makebox(0,0)[b]{$z\equiv y\:\wedge\: y\equiv x$}}
\put(240,63){\makebox(0,0)[b]{$y\equiv z\:\wedge\: x\equiv y$}}

\put(240,97){\makebox(0,0)[b]{$x\equiv y\:\wedge\: y\equiv z$}}

\put(193,17){\makebox(0,0)[r]{$t_{z,y,x}$}}
\put(193,49){\makebox(0,0)[r]{$s_{y,z}\wedge s_{x,y}$}}
\put(193,85){\makebox(0,0)[r]{$c_{x\equiv y,y\equiv z}$}}

\end{picture}
\end{center}

\noindent we can conclude that $G$ is indeed a functor.

To prove the faithfulness of $G$, i.e.\ coherence for $S_\equiv$,
we proceed in principle as for $S_\leq$ in Section~4. Now that we
have the equation ($ts$), we can permute freely
$t_{x,y,z}$-factors with $s_{x,z}$-factors, and eschew all the
complications we had with $M_\equiv$ in the preceding section.

\section{Preorder, equivalence and adjunction}

In this section we will show that assumptions concerning $\leq$
and $\equiv$ in categories of the preceding five sections amount
to assumptions about some adjoint situations. This matter is
related to matters considered in \cite{L70}.

Let $M^{-y}_\leq$ be the full subcategory of $M_\leq$ in whose
objects a particular variable $y$ does not occur, and let
$M^y_\leq$ be the full subcategory of $M_\leq$ whose objects are
of the form $y\leq u\wedge A$ for $y$ distinct from $u$ and not
occurring in $A$. Note that if the generating set $\cal V$ of
variables is infinite, then $M^{-y}_\leq$ and $M_\leq$ are
isomorphic categories.

For every variable $z$ distinct from $y$ there is a functor $F^z$
from $M^{-y}_\leq$ to $M^y_\leq$ defined as follows:

\begin{tabbing}

\hspace{12em}\=$F^zA\;$\=$=_{\df}\;\;$\=$y\leq z\wedge A$,
\\[1ex]
\>$F^zf$\>$=_{\df}$\>$\mj_{y\leq z}\wedge f$.

\end{tabbing}

\noindent Conversely, there is a functor $G^z$ from $M^y_\leq$ to
$M^{-y}_\leq$ where

\[
G^z(y\leq u\wedge A)=_{\df}\;z\leq u\wedge A,
\]

\noindent and $G^zf$ is obtained from the arrow term $f$ by
substituting $z$ for $y$, i.e.\ by uniformly replacing $y$ by $z$.
(The function $G^z$ on objects is also substitution of $z$ for
$y$, since $y$ is distinct from $u$ and does not occur in $A$.)

Let the arrow ${\gamma^z_A\!:A\vdash z\leq z\wedge A}$ of
$M^{-y}_\leq$, whose target is ${G^zF^zA}$, be defined by

\[
\gamma^z_A=_{\df}\;(r_z\wedge\mj_A)\cirk\sigma^\leftarrow_A,
\]

\noindent and let the arrow

\[
\varphi^z_{y\leq u\wedge A}\!:y\leq z\wedge(z\leq u\wedge A)\vdash
y\leq u\wedge A
\]

\noindent of $M^y_\leq$, whose source is ${F^zG^z(y\leq u\wedge
A)}$, be defined by

\[
\varphi^z_{y\leq u\wedge A}=_{\df}\;(t_{y,z,u}\wedge\mj_A)\cirk
b^\rightarrow_{y\leq z,z\leq u,A}.
\]

\noindent Then we can verify easily by appealing to coherence for
$M_\leq$ that the functor $F^z$ is left adjoint to $G^z$; in this
adjunction $\gamma^z$ is the unit natural transformation, and
$\varphi^z$ the counit natural transformation (see \cite{ML71},
Section IV.1).

The ``straightening of a sinuosity'' involved in the equations
$(rt\delta)$ and $(rt\sigma)$ (see the diagrams of Section~3)
indicated that we have such an adjunction (cf.\ \cite{D99},
Section 4.10, \cite{DP05}, Section 2.3, and \cite{D03}, Section
7).

The arrow $\gamma^z_A$ was defined in terms of $r_z$, but in
$M^{-y}_\leq$ we can define $r_z$ in terms of $\gamma^z$ as
follows:

\[
r_z=_{\df}\; \delta^\rightarrow_{z\leq z}\cirk \gamma^z_\top.
\]

\noindent Analogously, the arrow $\varphi^z_{y\leq u\wedge A}$ was
defined in terms of $t_{y,z,u}$, but we can define $t_{v,z,u}$ in
$M^{-y}_\leq$ in terms of $\varphi^z$ as follows. Note first that
in $M^y_\leq$ we can take

\[
t_{y,z,u}=_{\df}\;\delta^\rightarrow_{y\leq u}
\cirk\varphi^z_{y\leq u\wedge\top}\cirk(\mj_{y\leq
z}\wedge\delta^\leftarrow_{z\leq u})
\]

\noindent for $y$ different from $z$; in $M^{-y}_\leq$ we take

\[
t_{v,z,u}=_{\df}\;\delta^\rightarrow_{v\leq u} \cirk
G^v\varphi^z_{y\leq u\wedge\top}\cirk(\mj_{v\leq
z}\wedge\delta^\leftarrow_{z\leq u}).
\]

Then the specific equations of $M_\leq$ can be derived from the
assumption that we have the adjunction above between $M^{-y}_\leq$
and $M^y_\leq$, together with the equations

\begin{tabbing}

\centerline{$\gamma^z_A=((\delta^\rightarrow_{z\leq
z}\cirk\gamma^z_\top)\wedge\mj_A)\cirk\sigma^\leftarrow_A$,}
\\[2ex]
\centerline{$\varphi^z_{y\leq u\wedge
A}=((\delta^\rightarrow_{y\leq u}\cirk\varphi^z_{y\leq
u\wedge\top}\cirk(\mj_{y\leq z}\wedge\delta^\leftarrow_{z\leq
u}))\wedge\mj_A)\cirk b^\rightarrow_{y\leq z,z\leq u,A}$.}

\end{tabbing}

\noindent These equations are obtained from the definition of
$\gamma^z_A$ in terms of $r^z$ and the definition of
$\varphi^z_{y\leq u\wedge A}$ in terms of $t_{y,z,u}$. They give a
definition of $\gamma^z_A$ for an arbitrary $A$ in terms of
$\gamma^z_\top$, and a definition of $\varphi^z_{y\leq u\wedge A}$
for an arbitrary $A$ in terms of $\varphi^z_{y\leq u\wedge\top}$.

Note that a category equivalent to $M^y_\leq$ is the full
subcategory of $M_\leq$ whose objects have a single occurrence of
$y$ as the leftmost variable.

The adjunction we had above between $M^{-y}_\leq$ and $M^y_\leq$
is obtained also when $M_\leq$ is replaced by $S_\leq$, $M_\equiv$
and $S_\equiv$. With $M_\equiv$ we can take instead of the
category $M^y_\equiv$, defined like $M^y_\leq$, the equivalent
full subcategory $M^{y\ast}_\equiv$ of $M_\equiv$ whose objects
are those of the form $y\equiv u\wedge A$ or $u\equiv y\wedge A$
for $y$ distinct from $u$ and not occurring in $A$. We obtain an
adjunction between $M^{-y}_\equiv$, defined like $M^{-y}_\leq$
starting from $M_\equiv$, and $M^{y\ast}_\equiv$. 



With the categories $S^{-y}_\equiv$ and $S^{y\ast}_\equiv$,
defined analogously starting from $S_\equiv$, an analogous
adjunction obtains. The category $S^{y\ast}_\equiv$ is equivalent
to the full subcategory of $S_\equiv$ in whose objects $y$ occurs
exactly once.

\section{The coherence of $\dot S_\leq$}

Let us now suppose that terms are built with the help of a symbol
$\cdot$ that stands for a binary operation. We suppose, namely,
that terms are not only variables, but for $t_1$ and $t_2$ terms
we have that $t_1\cdot t_2$ is a term. We use $t$, $s$,
$r,\ldots,$ also with indices, for terms. To define atomic
formulae we suppose that if $t_1$ and $t_2$ are terms, then
$t_1\leq t_2$ is an atomic formula, as well as $\top$. Formulae
are defined otherwise as in Section~2.

The objects of the category $\dot S_\leq$ are these new formulae;
otherwise, $\dot S_\leq$ is defined like $S_\leq$ in Section~3
with the additional primitive arrow terms

\[
a_{t_1,t_2,t_3,t_4}\!\!:t_1\leq t_2\wedge t_3\leq t_4\vdash
t_1\cdot t_3\leq t_2\cdot t_4
\]

\noindent for all terms $t_1$, $t_2$, $t_3$ and $t_4$, which are
subject to the following additional equations:

\begin{tabbing}

\hspace{1em} \= $(ra)$\hspace{1em}\= $a_{t,t,s,s}\cirk(r_t\wedge
r_s)\cirk\delta^\leftarrow_\top=r_{t\cdot s}$,
\\[3ex]
for $c^m_{A,B,C,D}=_{\df}\; b^\rightarrow_{A,C,B\wedge
D}\cirk(\mj_A\wedge(b^\leftarrow_{C,B,D}\cirk(c_{B,C}\wedge\mj_D)\cirk
b^\rightarrow_{B,C,D}))\cirk b^\leftarrow_{A,B,C\wedge D}\!\!:$
\\[.5ex]
\` $(A\wedge B)\wedge(C\wedge D)\vdash(A\wedge C)\wedge(B\wedge
D)$,
\\[.5ex]
\> $(ta)$ \> $a_{t_1,r_1,t_2,r_2}\cirk(t_{t_1,s_1,r_1}\wedge
t_{t_2,s_2,r_2})=$
\\[.25ex]
\`$t_{t_1\cdot t_2,s_1\cdot s_2,r_1\cdot
r_2}\cirk(a_{t_1,s_1,t_2,s_2}\wedge a_{s_1,r_1,s_2,r_2})\cirk
c^m_{t_1\leq s_1,s_1\leq r_1,t_2\leq s_2,s_2\leq r_2}$.

\end{tabbing}

The equation $(ra)$ is parallel to the equation ($\wedge\:$1) of
Section~2, and $(ta)$ is parallel in the same manner to the
equation ($\wedge\:$2). In another manner, the equation $(ta)$ is
analogous to the equation (\emph{\v b\v c\v w}) of \cite{DP04}
(see the List of Equations).

We define the functor $G$ from $\dot S_\leq$ to the category
\emph{Br} by extending the definition of $G$ from $S_\leq$ to
\emph{Br} with a clause corresponding to the following diagram:

\begin{center}
\begin{picture}(100,40)
\put(30,9){\line(0,1){20}} \put(95,9){\line(0,1){20}}
\put(52,9){\line(1,1){20}} \put(72,9){\line(-1,1){20}}

\put(65,7){\makebox(0,0)[t]{$t_1\,\cdot\, t_3\:\leq\: t_2\,\cdot\,
t_4$}} \put(65,33){\makebox(0,0)[b]{$t_1\leq t_2\, \wedge\,
t_3\leq t_4$}}
\put(20,20){\makebox(0,0)[r]{$a_{t_1,t_2,t_3,t_4}$}}
\end{picture}
\end{center}

\noindent (where each line stands for a family of parallel lines;
cf.\ Section~3). To prove coherence for $\dot S_\leq$ we use
essentially the equation

\begin{tabbing}

\hspace{1em}\= $a_{t_1,s_1,t_2,s_2}\cirk(t_{t_1,r,s_1}\wedge
\mj_{t_2\leq s_2})=$
\\[.25ex]
\hspace{2em}\= $t_{t_1\cdot t_2,r\cdot s_2,s_1\cdot
s_2}\cirk(a_{t_1,r,t_2,s_2}\wedge a_{r,s_1,s_2,s_2})\cirk
c^m_{t_1\leq r,r\leq s_1,t_2\leq s_2,s_2\leq s_2}\cirk$
\\[.25ex]
\`$\cirk(\mj_{t_1\leq r\wedge r\leq s_1}\wedge((\mj_{t_2\leq
s_2}\wedge r_{s_2})\cirk\delta^\leftarrow_{t_2\leq s_2}))$,

\end{tabbing}

\noindent and another analogous equation with
${t_{t_1,r,s_1}\wedge \mj_{t_2\leq s_2}}$ on the left-hand side
replaced by ${\mj_{t_1\leq s_1}\wedge t_{t_2,r,s_2}}$; we proceed
otherwise using the ideas indicated in Section~3. This proof of
coherence is parallel to the proof of the Development Lemma of
Section~3; it formalizes the proof of this lemma on a different
level.

When we take the category $\dot S_\equiv$ obtained from $S_\equiv$
as $\dot S_\leq$ was obtained from $S_\leq$, with the additional
equation

\[
s_{t_1\cdot t_2,s_1\cdot s_2}\cirk
a_{t_1,s_1,t_2,s_2}=a_{s_1,t_1,s_2,t_2}\cirk(s_{t_1,s_1}\wedge
s_{t_2,s_2}),
\]

\noindent which is parallel to ${(f\wedge g)^{-1}=f^{-1}\wedge
g^{-1}}$, we can prove coherence analogously. In $\dot S_\equiv$
the relation corresponding to $\equiv$ is a congruence relation,
which, due to the scarcity of the means of expression of $\dot
S_\equiv$, is an equality relation.

To obtain coherence results for various categories extending $\dot
S_\equiv$, which would formalize fragments of the equational
theory of semigroups, or of the equational theory of monoids,
commutative or not, we would need additional arrows analogous to
the arrows $b^\rightarrow_{A,B,C}$, $b^\leftarrow_{A,B,C}$,
$\delta^\rightarrow_A$, $\delta^\leftarrow_A$,
$\sigma^\rightarrow_A$, $\sigma^\leftarrow_A$ and $c_{A,B}$. For
example, an arrow of type ${\top\vdash t_1\cdot(t_2\cdot
t_3)\equiv (t_1\cdot t_2)\cdot t_3}$ would correspond to
$b^\rightarrow$. The additional equations for these new arrows
would be parallel to the equations of monoidal or symmetric
monoidal categories.

\end{document}